\documentclass[12pt]{article}
\usepackage{amsfonts}
\usepackage{amsmath}
\usepackage{amssymb}
\usepackage{graphicx}

\usepackage{amsfonts}
\usepackage{amsmath}
\newtheorem{theorem}{Theorem}[section]

\newtheorem{corollary}[theorem]{Corollary}

\def\IH{{\mathbb H}} 
 
\def\IR{{\mathbb R}}

\def\IC{{\mathbb C}} 
\def\ID{{\mathbb D}}

\def\qed{ $\Box$}
\def\em{\it}

\title{Quasihyperbolic Metric Universal Covers}
\author{David A.\ Herron and Gaven J.\ Martin \thanks{ This research was initiated while both authors were supported by the Polish Academy of Sciences at Bedlewo.  We thank them for their hospitality. GJM is supported in part by a grant from the NZ Marsden Fund.\newline
\newline
{\bf  Mathematics subject classification 2010,}{Primary: 30F45, 30L99; Secondary:  51F99, 30C62, 20F65, 53A30, 53C23} \newline 
{\bf  Keywords,} {quasihyperbolic metric, universal cover, CAT(0), Hadamard}
}
\\ \\
\footnotesize{{\em Dedicated to Pekka Koskela on the ocassion of his $60^{\rm th}$ birthday}.}}
\date{}
\begin{document}

\maketitle

\begin{abstract}
We present a simple analytical proof that the natural metric universal cover of a quasihyperbolic planar domain is a complete Hadamard metric space.
\end{abstract}

\section{Introduction}

In an early paper \cite{Kos98}, Pekka Koskela wrote about old and new results regarding the quasihyperbolic metric.  Since then this metric has been a key player in much of Pekka's work, either as a fundamental geometric tool or as the focus of his research.  We anticipate that Pekka will find the following of interest.

Throughout this article $\Omega$ denotes a quasihyperbolic plane domain in the complex plane $\IC$.  Thus $\Omega$ is open and connected and $\Omega\neq \IC$.
Each such $\Omega$ carries a \emph{quasihyperbolic metric} $\delta^{-1}ds=\delta_\Omega^{-1}ds$ whose length distance $k=k_\Omega$ is called \emph{quasihyperbolic distance} in $\Omega$; here $\delta(z)=\delta_\Omega(z):={\rm dist}\,(z,\partial \Omega )$ is the Euclidean distance from $z$ to the boundary of $\Omega$.  See \ref{ss:qhyp m&d} below for more details.

\medskip

Quasihyperbolic distance in Euclidean domains was introduced by Gehring and Palka in \cite{GP76} and has since become an indispensable tool in the program of metric space analysis.
This paper is a companion to \cite{BH21}, \cite{Her21b}, \cite{HL21} where our program is to compare and contrast hyperbolic distance versus quasihyperbolic distance in planar domains.

A powerful tool available to hyperbolic geometers is the well-known fact that for any hyperbolic plane domain $\Omega$, there is a  holomorphic universal covering $(\IH^2,\tilde{h})\to(\Omega,h)$ which is a local isometry, where $(\IH^2,\tilde{h})$ is the standard hyperbolic plane and $h$ is hyperbolic distance in $\Omega$.  In \cite[Theorem~A]{Her21b} the first author established an analogous result for quasihyperbolic geometry.

Let $\tilde{\Omega}\stackrel{\Phi}{\to}\Omega$ be a  holomorphic universal cover of $\Omega$; topologically, $\tilde{\Omega}$ is a simply connected surface. 
When $\Omega$ is the once punctured plane $\Omega=\IC\setminus\{z_0\}$, $\tilde{\Omega}=\IC$ and $\Phi(z)=e^z+z_0$; otherwise, we can take $\tilde{\Omega}$ to be the unit disk. There is a unique length distance $\tilde{k}$ on $\tilde{\Omega}$ such that $(\tilde{\Omega},\tilde{k})\stackrel{\Phi}{\to}(\Omega,k)$ is a local isometry (see \cite[Prop.~3.25, p.42]{BH99} or \cite[p.80]{BBI01}) and $\tilde{k}$ is given by
\[
  \tilde{k}(a,b):=\inf_{\gamma} \ell_{\tilde{k}}(\gamma) \quad\text{where}\quad
  \ell_{\tilde{k}}(\gamma) := \ell_k(\Phi\circ\gamma)
\]
and the infimum is taken over all rectifiable paths $\gamma$ in $\tilde{\Omega}$ with endpoints $a, b$.

We give a simple proof of the following theorem that, essentially, relies solely on the fact that $-\log\delta$ is always subharmonic in $\Omega$.

\begin{theorem}  \label{TT:main}
For any quasihyperbolic plane domain $\Omega$, $(\tilde{\Omega},\tilde{k})$ is a Hadamard space.
\end{theorem}                   
\noindent
Recall that a \emph{Hadamard space} is a complete geodesic CAT(0) metric space.

 \medskip

As described in \cite[Corollary B]{Her21a}, there are some immediate consequences of the above.

\begin{corollary}  \label{CC:sc}  Let $\Omega$ be a simply connected quasihyperbolic plane domain.  Then

\begin{enumerate} 
  \item   $\Omega$, $(\Omega,k)$ is Hadamard.
  \item  For each pair of points $a,b$ in $\Omega$, there is a unique quasihyperbolic geodesic with endpoints $a,b$.
  \item  For each pair of points $a,b$ in any quasihyperbolic plane domain, each homotopy class of paths in $\Omega$ with endpoints $a,b$ contains a unique quasihyperbolic geodesic.
\end{enumerate}
\end{corollary}                  

In \cite[Theorem~C]{Her21a} we proved that the injectivity radius of $(\tilde{\Omega},\tilde{k})\stackrel{\Phi}{\to}(\Omega,k)$ is at least $\pi$.  Here, comparison theorems give quick bounds,  but achieving the sharp bound requires more effort,  so we do not attempt it here.  We mention that already in \cite[Corollary 3.6, p.44]{MO86} Martin and Osgood proved that quasihyperbolic plane domains have non-positive generalized quasihyperbolic Gaussian curvature, and therefore the above Theorem is not so surprising. 

\section{Preliminaries}  
%
\subsection{General Information}  

We view the Euclidean plane as the complex number field $\IC$.  Everywhere $\Omega$ is a {\em quasihyperbolic plane domain}. 
The open unit disk is $\ID:=\{z\in\IC:|z|<1\}$ and $\IC^*:=\IC\setminus\{0\}$.
 The quantity $\delta(z)=\delta_\Omega(z):={\rm dist}\,(z,\partial\Omega)$ is the Euclidean distance from $z\in\IC$ to the boundary of $\Omega$, and $1/\delta$ is the scaling factor (that is the metric-density) for the so-called \emph{quasihyperbolic metric} $\delta^{-1}ds$ on $\Omega\subset\IC$; see \S \ref{ss:qhyp m&d}.

\subsection{Conformal Metrics}      \label{s:cfml metrics} %
Each positive continuous function $\rho:\Omega\to(0,+\infty)$ induces a length distance $d_\rho$ on $\Omega$ defined by
\[
  d_\rho(a,b):=\inf_{\gamma } \ell_\rho(\gamma) \quad\text{where}\quad \ell_\rho(\gamma):= \int_\gamma \rho\, ds
\]
and where the infimum is taken over all rectifiable paths $\gamma$ in $\Omega$ that join the points $a,b$.  We describe this by calling  $\rho\,ds=\rho(z)|dz|$ a \emph{conformal metric} on $\Omega$.  The geodesic regularity theorem of \cite{Mar} shows that as soon as $\rho$ is continuous and the metric is complete the infimum is attained by a curve with Lipschitz continuous first derivatives - a $\rho$ geodesic.

When $\rho$ is sufficiently differentiable,  say $C^2(\Omega)$, the \emph{Gaussian curvature} of $(\Omega,d_\rho)$ is given by
\[
{\bf K}_\rho:=-\rho^{-2}\Delta\log\rho\,.
\] 
In the plane,  this curvature can also be defined for continuous metric densities through an integral formula for the Laplacian -  though it is not always finite,  see \cite[\S 3]{MO86}. This idea was first observed by Heins \cite{Heins} who proved a version of the Schwarz lemma when the weak curvature was $\leq -1$.

\subsubsection{The QuasiHyperbolic Metric \& Distance} 	\label{ss:qhyp m&d} %
The \emph{quasihyperbolic distance} $k=k_\Omega$ in $\Omega$ is the length distance $k_\Omega:=d_{\delta^{-1}}$ induced by the quasihyperbolic metric $\delta^{-1}ds$ on $\Omega$ where $\delta(z)=\delta_\Omega(z):={\rm dist}\,(z,\partial\Omega)$ is the Euclidean distance from $z$ to the boundary of $\Omega$.
Thanks to the Hopf-Rinow Theorem (see \cite[p.35]{BH99}, \cite[p.51]{BBI01}, \cite[p.62]{Pap05}) we know that $(\Omega,k)$ is a geodesic metric space. 
Indeed, any rectifiably connected non-complete locally compact metric space admits a quasihyperbolic metric that is geodesic.  See also \cite[Lemma~1]{GO79} for Euclidean domains. Further,  the  basic estimates for quasihyperbolic distance were established by Gehring and Palka \cite[2.1]{GP76}: For all $a,b\in\Omega$,
\[  k(a,b) \geq \log\left(1+\frac{l(a,b)}{\delta(a)\wedge \delta(b)}\right)  \geq  \log\left(1+\frac{|a-b|}{\delta(a)\wedge \delta(b)}\right) \ge\left|\log\frac{\delta(a)}{\delta(b)}\right|\, \]
 here $l(a,b)=l_\Omega(a,b)$ is the (intrinsic) Euclidean length distance between $a$ and $b$ in $\Omega$.  The first inequality above is a special case of the more general (and easily proven) inequality
\[   \ell_k (\gamma) \ge \log\bigl(1+\ell(\gamma)/{\rm dist}\,(|\gamma|,\partial\Omega) \bigr)  \]
 which holds for any rectifiable path $\gamma$ in $\Omega$.  Using these estimates, we deduce that
 if $k(a,b)\leq 1$ or $  |a-b| \leq\frac{1}{2}\delta(a)$\,, then              \[      \frac{1}{2} k(a,b)\leq\frac{|a-b|}{\delta(a)}\leq 2 k(a,b) . \]
These estimates reveal that $(\Omega,k)$ is complete.

In another direction, if $\gamma:[0,L]\to\Omega$ is a quasihyperbolic geodesic in $\Omega$ parameterized  with respect to  Euclidean arclength, then $\gamma$ is differentiable (even at its endpoints) and for all $s,t\in[0,L]$, 
 \[  |\gamma'(s)-\gamma'(t)| \le 2\, \frac{|\gamma(s)-\gamma(t)|}{\delta\bigl(\gamma(s)\bigr)}\, \leq  \frac{2e^\frac{1}{2}}{\delta(a)} |s-t|\,.\]
  Uniform estimates such as this,  along with comparison principles,  establish the convergence of geodesics in $C^{1,1}$ for locally uniformly convergent sequences of complete continuous metric densities,  in particular for our smoothing approach via mollifiers given below.  We leave the reader to explore this.

As a well known example,  well note that $(\IC^*,k_*)$ is (isometric to) the Euclidean cylinder $\mathbb{S}^1\times\IR^1$ with its Euclidean length distance inherited from its standard embedding into $\IR^3$.  
One way to realize this is via the holomorphic covering $\exp:\IC\to\IC^*$ which pulls back the quasihyperbolic metric $\delta_*^{-1}ds$ on $\IC^*$ to the Euclidean metric on $\IC$, as explained in \cite{MO86}.  In particular, quasihyperbolic geodesics in $\IC^*$ are logarithmic spirals, for all $a,b\in\IC^*$,
\[
  k_*(a,b)=\bigl|\log(b/a)\bigr|=\bigl|\log|b/a|+i\arg(b/a)\bigr|\,,
\]
and $\exp:(\IC,|\cdot|)\to(\IC^*,k_*)$ is a metric universal cover.  Thus our Theorem holds when $\Omega$ is a once punctured plane.

\subsection{CAT(0) Metric Spaces}      \label{s:CAT0} %
Here our terminology and notation conforms exactly with that in \cite{BH99} and we refer the reader to this delightful trove of geometric information about non-positive curvature, and also see \cite{BBI01}.  We recall a few fundamental concepts, mostly  copied directly from \cite{BH99}.  Throughout this subsection, $X$ is a geodesic metric space; for example, $X$ could be a quasihyperbolic plane domain with its quasihyperbolic distance, or a closed rectifiably connected plane set with its intrinsic length distance.

\subsubsection{Geodesic and Comparison Triangles}  \label{ss:triangles} 
A \emph{geodesic triangle} $\Delta$ in $X$ consists of three points in $X$, say $a,b,c\in X$, called the \emph{vertices of $\Delta$} and three geodesics, say $\alpha:a\curvearrowright b, \beta:b\curvearrowright c, \gamma:c\curvearrowright a$ (that we may write as $[a,b], [b,c], [c,a]$) called the \emph{sides of $\Delta$}.  We use the notation
\[
  \Delta=\Delta(\alpha,\beta,\gamma) \quad\text{or}\quad \Delta=[a,b,c]:=[a,b]\star[b,c]\star[c,a] \quad\text{or}\quad \Delta=\Delta(a,b,c)
\]
depending on the context and the need for accuracy.

A Euclidean triangle $\bar\Delta=\Delta(\bar{a},\bar{b},\bar{c})$ in $\IC$ is a \emph{comparison triangle} for $\Delta=\Delta(a,b,c)$ provided $|a-b|=|\bar{a}-\bar{b}|,|b-c|=|\bar{b}-\bar{c}|,|c-a|=|\bar{c}-\bar{a}|$.  We also write $\bar\Delta=\bar\Delta(a,b,c)$ when a specific choice of $\bar{a},\bar{b},\bar{c}$ is not required.  A point $\bar{x}\in[\bar{a},\bar{b}]$ is a \emph{comparison point} for $x\in[a,b]$ when $|x-a|=|\bar{x}-\bar{a}|$.

\subsubsection{CAT(0) Definition}  \label{ss:CAT0} 

A geodesic triangle $\Delta$ in $X$ satisfies the \emph{CAT(0) distance inequality} if and only if the distance between any two points of $\Delta$ is not larger than the Euclidean distance between the corresponding comparison points; that is,
\[
  \forall\; x,y\in\Delta\;\text{and corresponding comparison points}\;\bar{x},\bar{y}\in\bar\Delta\;, \quad |x-y| \le |\bar{x}-\bar{y}|\,.
\]
A geodesic metric space is \emph{CAT(0)} if and only if each of its geodesic triangles satisfies the CAT(0) distance inequality.

A complete CAT(0) metric space is called a \emph{Hadamard space}.  A geodesic metric space $X$ has \emph{non-positive curvature} if and only if it is locally CAT(0), meaning that for each point $a\in X$ there is an $r>0$ (that can depend on $a$) such that the metric ball $B(a;r)$ (endowed with the distance from $X$) is CAT(0).

Each sufficiently smooth Riemannian manifold has non-positive curvature if and only if all of its sectional curvatures are non-positive; see \cite[Theorem~1A.6, p.173]{BH99}.  In particular, if $\rho\,ds$ is a smooth conformal metric on $\Omega$ with $\bf K_\rho\le0$, then $(\Omega,d_\rho)$ has non-positive curvature.

\subsection{Smoothing}      \label{ss:smooth} %
As usual, we start with a $C^\infty(\Omega)$ smooth $\eta:\IC\to\IR$ with $\eta\ge0$, $\eta(z)=\eta(|z|)$, the support of $\eta$ lies in $\ID$, and $\int_\IC\eta=1$.  For each $\varepsilon >0$ we set $\eta_\varepsilon (z):=\varepsilon ^{-2}\eta(z/\varepsilon )$.  The \emph{regularization} (or \emph{mollification}) of an $L^{1}_{loc}(\Omega)$ function $u:\Omega\to\IR$ are the convolutions $u_\varepsilon :=u\ast\eta_\varepsilon $, so
\[
  u_\varepsilon (z) := \int_\IC u(w)\eta_\varepsilon (z-w)\, dA(w)\,,
\]
which are defined in $\Omega_\varepsilon :=\{z\in\Omega:\delta(z)>\varepsilon \}$.  It is well known that $u_\varepsilon \in C^\infty(\Omega_\varepsilon )$ and $u_\varepsilon \to u$ as $\varepsilon \to0^+$ where this convergence is: pointwise at each Lebesgue point of $u$, locally uniformly in $\Omega$ if $u$ is continuous in $\Omega$, and in $L^p_{loc}(\Omega)$ if $u\in L^p_{loc}(\Omega)$.  Moreover, if $u$ is subharmonic in $\Omega$, then so is each $u_\varepsilon $.  See for example \cite[Proposition I.15, p.235]{GL86} or \cite[Theorem~2.7.2, p.49]{Ran95}.

\section{Proof of Theorem}  \label{S:PfThm} 

Let $\Omega\subsetneq\IC$ be a planar domain.  The main idea is to approximate $(\Omega,k)$ by metric spaces $(\Omega_\varepsilon ,d_\varepsilon )$ that all have non-positive curvature.  Then a limit argument, exactly like that used in the proof of \cite[Theorem~A]{Her21a}, gives the asserted conclusion.

We start with the well known fact that $u(z):=-\log\delta(z)$ is subharmonic in $\Omega$.  This is easy to see.  Evidently, $u$ is continuous.  For each fixed $\zeta\in\partial\Omega$, $z\mapsto-\log|z-\zeta|$ is harmonic in $\IC\setminus\{\zeta\}\supset\Omega$ and so has the mean value property in $\Omega$.  It follows that $u$ has the submean value property in $\Omega$.

Let $u_\varepsilon :=u\ast\eta_\varepsilon $ be the regularization of $u$ as described at (\ref{ss:smooth}).  Thus $u_\varepsilon $ is defined, $C^\infty$ smooth, and subharmonic (so, $\Delta u_\varepsilon \geq 0$) in $\{z\in\Omega:\delta(z)>\varepsilon \}$.  Moreover,  $u_\varepsilon \to u$ as $\varepsilon \to0^+$ locally uniformly in $\Omega$.

We assume that $\Omega$ is not a once punctured plane and that the origin lies in $\Omega$.  Put $\varepsilon _n:=\delta(0)/n, u_n:=u_{\varepsilon _n}$, and let $\Omega_n$ be the component of $\{z\in\Omega:\delta(z)>\varepsilon _n\}$ that contains the origin.  Then $\{\Omega_n\}_{n=1}^{\infty}$ Carath\'eodory kernel converges to $\Omega$  with respect to the origin.

Next, let $\rho_n:=e^{u_n}$.  Then $\rho_n>0$ and $C^\infty$ in $\Omega_n$.  Let $d_n:=d_{\rho_n}$ be the length distance in $\Omega_n$ associated to the conformal metric $\rho_n\,ds$.  Since $u_n$ is subharmonic in $\Omega_n$, $\rho_n\,ds$ has Gaussian curvature
\[
{\bf K}_{\rho_n}=-\rho_n^{-2}\Delta\log\rho_n\le0\quad\text{in $\Omega_n$}\,.
\]
It follows that the metric spaces $(\Omega_n,d_n)$ all have non-positive curvature; see \cite[Theorem~1A.6, p.173]{Her21a}.

Let $ {\Phi}:\ID\to \Omega$ and ${\Phi_n}:\ID\to \Omega_n$ be holomorphic covering projections with $\Phi(0)=0=\Phi_n(0)$ and $\Phi'(0)>0, \Phi'_n(0)>0$.  Since $\{\Omega_n\}_{n=1}^{\infty}$ Carath\'eodory kernel converges to $\Omega$  with respect to the origin, a theorem of Hejhal's \cite[Theorem~1]{Hej74} 
(see also \cite[Corollary 5.3]{BM22}) 
asserts that $\Phi_n\to \Phi$, so also $\Phi'_n\to \Phi'$, locally uniformly in $\ID$.

Let $\tilde{k}, \tilde{d}_n$ be the $\Phi,\Phi_n$ lifts of the distances $k,d_n$ on $\Omega, \Omega_n$ respectively.  That is, $\tilde{k}$ and $\tilde{d}_n$ are the length distances on $\ID$ induced by the pull backs
\begin{gather*}
  \tilde\rho\,ds:=\Phi^*\bigl[\delta^{-1}ds] \quad\text{and}\quad  \tilde\rho_n\,ds:=\Phi_n^*\bigl[\rho_n\,ds]
  \intertext{of the metrics $\delta^{-1}ds$ and $\rho_n\,ds$ in $\Omega$ and $\Omega_n$ respectively.  Thus, for $\zeta\in\ID$,}
  \tilde\rho(\zeta)\,|d\zeta|=\frac{|\Phi'(\zeta)|}{\delta\bigl(\Phi(\zeta)\bigr)}\,|d\zeta|  \quad\text{and}\quad  \tilde\rho_n(\zeta)\,|d\zeta|=|\Phi_n'(\zeta)| \, \rho_n\bigl(\Phi_n(\zeta)\bigr)\,|d\zeta|
\end{gather*}
and $\Phi: (\ID,\tilde{k})\to(\Omega,k), \; \Phi_n(\ID,\tilde{d}_n)\to(\Omega_n,d_n)$ are metric universal coverings.

Note that as $(\Omega_n,d_n)$ has non-positive curvature, the Cartan-Hadamard Theorem \cite[Chapter~II.4, Theorem~4.1, p.193]{BH99} asserts that $(\ID,\tilde{d}_n)$ is CAT(0).

Using the locally uniform convergences of $\rho_n\to\delta^{-1}$ and $\Phi_n\to\Phi, \Phi'_n\to\Phi'$ (in $\Omega$ and $\ID$ respectively) we deduce that $\tilde\rho_n\,ds\to\tilde\rho\,ds$ locally uniformly in $\ID$. This implies pointed Gromov-Hausdorff convergence of $(\ID,\tilde{d}_n,0)$ to $(\ID,\tilde{k},0)$ (see the proof of \cite[Theorem~4.4]{HRS20}) which in turn says that $(\ID,\tilde{k})$ is a 4-point limit of $(\ID,\tilde{d}_n)$ and hence, as each $(\ID,\tilde{d}_n)$ is CAT(0), it follows that $(\ID,\tilde{k})$ is CAT(0); see \cite[Cor.~3.10, p.187; Theorem~3.9, p.186]{BH99}.  Finally, it is a routine matter to check that $(\ID,\tilde{k})$ is complete; for instance, see \cite[Exercise~3.4.8, p.~80]{BBI01}.
\hfill \qed

DAH,  Department of Mathematical Sciences, University of Cincinnati, OH 45221-0025, USA \\
email: David.Herron@UC.edu

\medskip

GJM, Institute for Advanced Study, Massey University, Auckland, New Zealand \\
email: G.J.Martin@Massey.ac.NZ

\end{document}